\documentclass[12pt,a4paper]{article} %[pwi.tex 05.06.22  <= 02.09.21]
\usepackage{amsfonts,amssymb} 
\usepackage[cp1251]{inputenc}

\usepackage[margin=25mm, a4paper]{geometry}  % journal quality

\usepackage{tikz}  \def\bwait#1#2{#2}

\def\n{\noindent}  \def\?#1{}
\def\IZ{{\mathbb{Z}}}  \def\IR{{\mathbb{R}}}
\def\IQ{{\mathbb{Q}}}   \def\cB{{\cal B}} 
        
   \def\Tor{{\mathbb{T}}}     
  \def\Clos{{\rm Clos}}  \def\Per{{\rm Per}}  \def\R{{\cal R}}

   \def\phi{\varphi} 
\def\ep{\varepsilon}  \def\t{\tilde} \def\bdelete#1{}
  \def\dist{{\rm dist}}
\def\mod1{\,({\rm mod\ } 1)\,}

\def\toas#1{\stackrel{#1}{\longrightarrow}}

\def\proof{\smallskip \noindent {\bf Proof. \ }}       %start of proof
\def\blanksquare{\,\,\,$\sqcup\!\!\!\!\sqcap$}         %blank  square
\def\qed{\hfill\blanksquare\linebreak\smallskip\par}   %end of proof
\def\hence{\Longrightarrow}                                 %==>

\def\thname{Theorem}  \def\lmname{Lemma}    \def\prname{Proposition}
\def\dfname{Definition}  \def\crname{Corollary}  \def\rmname{Remark}
\def\exname{Example}  

\newtheorem{theorem}{\thname}[section]   %Numbering: Theorem--Other section
\newtheorem{lemma}{\lmname}[section]     %{lemma}[theorem]{Lemma}   subsection
\newtheorem{proposition}[lemma]{\prname} %lemma
\newtheorem{example}{\exname}[section]

\newtheorem{dftn}{\dfname}[section]
\newenvironment{definition}{\begin{dftn}\rm}{\end{dftn}} %section
\def\bdef#1{\begin{definition} #1 \end{definition}}
\newtheorem{rmrk}[lemma]{\rmname}
\makeatletter\def\fps@figure{htbp}\makeatother %figure pos: tbp - standard

\begin{document}

\title{Invariant measures of torus piecewise isometries}
\author{Michael Blank\thanks{
        Institute for Information Transmission Problems RAS
        (Kharkevich Institute);}
        \thanks{National Research University ``Higher School of Economics'';
        e-mail: blank@iitp.ru}
       }
\date{September 2, 2021} %\today}
\maketitle

\begin{abstract}%
We study measure-theoretical aspects of torus piecewise isometries. 
Not much is known about this type of dynamical systems, except for 
the special case of one-dimensional interval exchange mappings. The 
last case is fundamentally different from the general situation 
in the presence of an invariant measure (Lebesgue measure), which 
helps a lot in the analysis. Due to the absence of good methods of 
analysis of general systems with discontinuities, even the existence of 
invariant measures of the torus piecewise isometries was an open question. 
We establish sufficient conditions for the existence/absence of invariant 
measures for this class of systems. Technically, our results are based 
on the approximation of the maps under study by weakly periodic ones. 
\end{abstract}% 

{\small\n
}

%%%%%%%%%%%%%%%%%%%%%%%%%%%%%%%%
\section{Introduction}\label{s:Intro}
By now, we have learned reasonably well how to study hyperbolic 
(locally expanding/contracting or both) chaotic dynamical systems, 
thanks to a large extent to the development of the so called 
operator approach. Contrary to this not much is known about 
piecewise isometries, except for a special case of one-dimensional 
interval exchange transformations (IET) and a few similar very 
exceptional multidimensional systems. It is worth noting that IETs are 
fundamentally different from the general situation in the clear presence 
of an invariant measure (Lebesgue measure), which helps a lot in the analysis. 

While the IET represent mainly pure theoretical constructions, the general 
piecewise isometries appear naturally in various applications like 
contemporary methods of machine learning, some piecewise smooth 
physical models (in particular Fermi-Ulam models), etc.  (see, e.g. 
\cite{De} for applications in engineering, and further references about 
other possibilities therein). Indeed, local translations of an IET disrupt 
its structure completely, while preserve the class of general piecewise 
isometries. Note also that the maps from this class in general are 
non-invertible, which adds possibilities for applications for ``real life'' 
modeling, but also additional problems for their analysis. 

In general from a measure-theoretical point of view one of the first steps 
in the analysis of a dynamical system is the study of its invariant measures. 
In some cases one can easily find a ``good'' invariant measure (e.g., for IETs), 
or to prove its existence if the map, defining the dynamical system, is 
continuous or satisfies some special properties (e.g., piecewise expanding). 
For systems with singularities (e.g., discontinuities) even the question about 
the existence of an invariant measure might become a difficult problem. 
The class of piecewise isometries represents an example of the latter sort. 
The point is that the methods (e.g., operator approach or symbolic dynamics) 
well established for other types of chaotic dynamical systems do not work 
in this setting and one needs to look for new approaches, one of which will 
be proposed in the present paper.

Let us start with a simple example of an isometry in $\IR^2$ defined by the map 
$\t{T}x := \R{x} + b$, where $x,b\in\IR^2$ are vectors and $\R$ is a rotation matrix. 
Applying the same map to the 2-dimensional torus $X$ we get $Tx:=\t{T}x \mod1$, 
which is no longer a global isometry. Moreover, the map $T$ is not a bijection and 
has both lines of discontinuities and regions of points having no preimages. On the 
other hand, one can easily make a partition of $X$ into a finite number of regions, 
such that the restriction of $T$ to each region is an isometry (see Section~\ref{s:affine} 
for details). This is a basic model for our study. 
Some nontrivial invariant sets for such a map are demonstrated at Fig.~\ref{f:linear}. 
We are interested in asymptotic properties of the dynamics of piecewise isometric 
torus maps like the one shown in this figure.

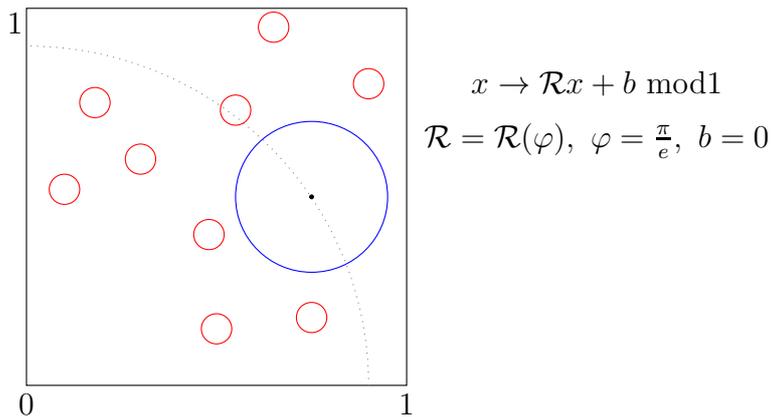
\begin{figure}
\begin{center}
\begin{tikzpicture}[scale=0.5]
       \draw (0,0) rectangle (10,10); 
       \draw[gray, dotted] (9,0) arc [radius=9, start angle=0, end angle=90];
       \bwait{3}{\draw[black, fill] (7.5,5) circle(.05);}
       \bwait{2}{\draw[blue] (7.5,5) circle(2);}
       \bwait{4}{
       \draw[red] (1,5.2) circle(.4); \draw[red] (1.8,7.5) circle(.4); \draw[red] (3,6) circle(.4);
       \draw[red] (4.8,4) circle(.4); \draw[red] (5,1.5) circle(.4); \draw[red] (7.5,1.8) circle(.4);
       \draw[red] (5.5,7.3) circle(.4); \draw[red] (6.5,9.5) circle(.4); \draw[red] (9,8) circle(.4);}
%       \draw [thick] (0,1) .. controls (2,3) and (3,4) .. (4.5,10);
%       \draw [thick] (10,1) .. controls (8,3) and (7,4) .. (5.5,10);
%      \draw [->] (5,0) -- (5,10); \draw [->] (5,0) -- (5,5);
       %
       \bwait{1}{
       \node at (0,-0.5) {$0$};  \node at (10,-0.5) {$1$};  \node at (-0.3,9.6) {$1$};}
       \bwait{1}{{\node at (15,8) {$x \to \R x + b$ mod1}; 
                              \node at (15,6.5) {$\R=\R(\phi), ~\phi=\frac\pi{e}, ~b=0$};}}
\end{tikzpicture}\end{center}
\caption{Invariant sets of the torus rotation by the angle $\pi/e$}\end{figure}
\label{f:linear} 

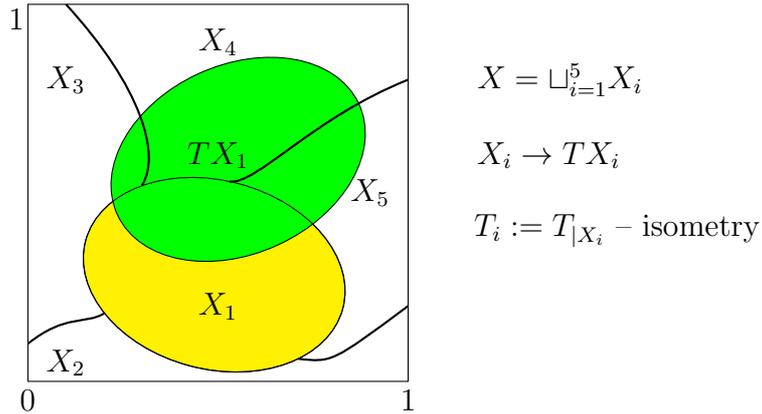
\begin{figure} 
\begin{center}
\begin{tikzpicture}[scale=0.5]
      \draw (0,0) rectangle (10,10); 
      \node at (0,-0.5) {$0$};  \node at (10,-0.5) {$1$};  \node at (-0.3,9.6) {$1$};
      \draw[fill=yellow,rotate=-15](4,4) ellipse (3.5cm and 2.5cm);
      \bwait{2}{
      \draw[fill=green,rotate=+25](6+1.5,4-1) ellipse (3.5cm and 2.5cm);
      \draw[rotate=-15](4,4) ellipse (3.5cm and 2.5cm); \node at (5,6){$TX_1$};}
      \draw [thick]  (0,1) .. controls (1,1.8) and (1.5,1.5) .. (2,1.8);  %Bezier
      \draw [thick]  (1,10) .. controls (3,7.8) and (3.5,6) .. (3,5.2);  %Bezier
      \draw [thick]  (5.3,5.3) .. controls (6,5.2) and (7.5,7) .. (10,8);  %Bezier
      \draw [thick]  (7.1,.6) .. controls (8,.5) and (8,.5) .. (10,2);  %Bezier
      \node at (5,2){$X_1$};    %\node[blue]
      \node at (1,.5){$X_2$}; \node at (1,8){$X_3$}; 
      \node at (5,9){$X_4$};  \node at (9,5){$X_5$}; 
      \node at (14,8){$X=\sqcup_{i=1}^5X_i$};   \node at (13.7,6){$X_i \to TX_i$}; 
      \node at (15.5,4){$T_i:=T_{|X_i}$ -- isometry}; 
\end{tikzpicture}\end{center}
\caption{A special partition for a torus isometry. 
One of the elements of the partition and its image are marked in yellow and 
green colors correspondingly.}\label{f:t-isometry} 
\end{figure}

Recall that an isometry is a map which preserves distances between all points in the 
phase space. Thus a piecewise isometry (PWI) is a map $T$ from a space $X$ equipped with 
a special partition $\{X_i\}$ into itself such that the restriction of $T$ to each element of the partition 
preserves distances (i.e. is a true isometry). An example of a torus PWI is demonstrated at 
Fig.~\ref{f:t-isometry}. 

Naturally, the list of available isometries sensitively depends on the choice of a metric. 
In particular, under the discrete metric each bijection turns out to be an isometry. 
We will consider a natural family of translational invariant metrics 
$\t\rho_p(x,y):=(\sum_i|x_i - y_i|^p)^\frac1p,~p\ge1$, 
where $x_i,y_i$ are the coordinates of the points $x,y\in X\in\IR^d$.
The Euclidean metric corresponds to $\t\rho_2(\cdot,\cdot)$, while the uniform 
metric to $\t\rho_\infty(\cdot,\cdot)$. The corresponding family of torus metrics  
is defined as $\rho_p(x,y) := \min_{z\in\IZ^d}\t\rho_p(x-z,y)$.

The paper is organized as follows. In Section~\ref{s:history} we give 
a brief historical description of results related to piecewise isometries 
and motivations of the present research. Section~\ref{s:setup} is 
dedicated to basic definitions and the formulation of the main result -- Theorem~\ref{t:main}. 
Examples demonstrating the possibility of the absence of invariant measures 
for PWI and the focusing phenomenon (due to which this happens) will be 
described in Section~\ref{s:absence}. 
In Section~\ref{s:per} we introduce and study weakly periodic maps and 
semigroups. In Section~\ref{s:appr} we apply the approximation by 
weakly periodic maps to prove Theorem~\ref{t:main}. 
Finally in Section~\ref{s:affine} we study a seemingly simple example of a PWI: 
torus restriction of a plane rotation. Note that this situation is excluded 
from Theorem~\ref{t:main} by the extendability assumption.

\bigskip 

\n{\bf Acknowledgments.} The author is grateful to Gregory Olshansky and  
Leonid Polterovich for useful discussions.

\section{Historical remarks} \label{s:history}
We start this short review with the discussion of known results related to ``bijective'' 
PWI $T:X\to X$, when the union of the images under $T$ of the elements of the special 
partition $\{X_i\}$ covers the entire phase space $X$. Since in the multidimensional 
setting ($d>1$) this situation is rather exotic, we will discuss only the one-dimensional 
case, known in the literature as interval exchange transformations (IET). 
Quotes marks in the word ``bijective'' emphasize a problem with the boundary points 
of the elements of the partition $\Gamma:=\cup_i\partial{X_i}$. 
In fact, the authors either assume that either $TX_i\cap TX_j\in T(\Gamma)$ for $i\ne j$ 
(may intersect only on boundary points), or instead of the partition consider the 
covering $\{\Clos{X_i}\}$. The former case leads to an endomorphism,\footnote{Except for 
   orientation preserving PWI with a very special structure of the partition, consisting of 
   semi-closed intervals of type $X_i:=[a_i,a_{i+1})$.} 
while the latter to a multi-valued map, because the map is not uniquely defined at $\Gamma$. 
This is especially unpleasant because under the action of a multi-valued map a total 
measure of the phase space may grow, which cannot happen for a conventional 
single-valued map. Fortunately, in the case of IET for any version of the action on 
boundary points the Lebesgue measure is preserved  under dynamics. Therefore, 
since $\Gamma$ is the set of zero Lebesgue measure, one does not care much about 
these peculiarities. 

Analysis of IETs is a very popular field of research and the list of publications related to 
them is very long and still actively growing. Therefore we mention only some basic 
publications in a historical order and refer the reader to look for further reference 
therein. The first results and important constructions related to orientation preserving 
IET were described in V.I.~Oseledets (1966) \cite{Os}, M.~Keane (1975) \cite{Ke}, 
A.~Katok (1980) \cite{Ka}, W.~Veech (1982) \cite{Ve}, H.~Masur (1982) \cite{Ma}, 
M.~Viana (2006) \cite{Vi}, A.I.~Bufetov( 2006) \cite{Bu}. 
Starting already from the work of W.~Veech the main object of interest here was shifted 
from a pure dynamical point of view to connections with geometry, in particular with 
Teichmuller flows. Since the latter goes out of the scope of the present paper, we will not discuss this. 

In the present paper we are interested mainly in measure-theoretic points of view 
to dynamics, which is relatively simple and well understood in the case of IETs.

The situation is becoming rather different if a PWI is not ``bijective'' even in the 
above mentioned weak sense. This means that the images of the elements of 
the special partition might have arbitrary intersections. The main body of results 
in this direction are again about the one-dimensional orientation preserving case, 
known in the literature as interval translation maps (ITM). 
The most interesting among them in a historical order belongs to  
M.~Boshernitzan \& I.~Kornfeld (1995) \cite{BK}, 
J.~Buzzi \& P. Hubert (2004) \cite{BH}, H.~Bruin (2007,2012) \cite{Br}, 
S.~Marmi, P.~Moussa, J-C.~Yoccoz (2012) \cite{MMY}, D.~Volk (2014) \cite{Vo}.
A further generalization when a general cover is used instead of a special partition 
(which inevitably leads to a multi-valued map) was considered in 
A.~Skripchenko \& S.~Troubetzkoy (2015) \cite{ST}. In any case the problem 
with the boundary points already mentioned in the case of IET is becoming 
a serious obstacle here. Nevertheless from a topological point of view or 
in the discussion of attractors one can go around it in order to reduce the 
analysis of ITM to a better undrestood IET. 
A popular approach introduced in \cite{BK} for the case of the circle orientation 
preserving PWI is to classify them into finite and infinite types with respect to the 
number of iterations under which the images of the phase space are stabilized. 
Namely, the finite type is characterized by the property $T^{n+1}X= T^{n}X$ 
for some $n<\infty$, and otherwise we refer to an infinite type. 
Indeed, for the finite type case in the one-dimensional orientation preserving case 
one obtains an ITM in a finite number of iterations. 

A number of attempts to follow this idea in the multidimensional setting have been 
tried (see, e.g. \cite{Vo,Tr} and further references therein). Unfortunately, in 
distinction to the one-dimensional case, the multidimensional dynamics is much 
more complicated and cannot be reduced to some version of ITM, which has been 
demonstrated by A.~Goetz in \cite{Go}. In his example $X$ is an isosceles triangle 
divided by the special partition into two smaller isosceles triangles $X_1$ and $X_2$. 
Under the action of the piecewise isometry $T$ the smaller triangles are rotated 
one clockwise and another anticlockwise in such a way that their images intersect 
only along boundaries and their closures cover $X$ completely. Similarly to the 
discussion above, there is a problem how to define the action of $T$ on the 
boundary of the special partition, but on the level of numerics this is not so 
important. The aim of the author was to study numerically the structure of 
limit sets under dynamics, which turns out to be extremely complicated. 
Moreover the simulation indicates that this structure is self similar up to a  
very fine scales. This conjecture is still not proven on the mathematical level. 

From the most interesting for us measure-theoretic point of view, there are no 
results at all for $d>1$ and only a few in the one-dimensional setting, which 
we will discuss now. The main problem here is that there is no obvious candidate 
for an invariant measure, while the images of a good enough reference 
measure in general might concentrate at boundary points, which complicates 
their analysis.  In \cite{BH} it has been shown that in the absence of periodic 
points the number of ergodic invariant measures of an orientation preserving 
circle PWI cannot be larger than the number of elements of the special partition. 
It is somewhat strange that having this answer the authors even did not 
consider a question whether at least one invariant measure is present. 
Observe that PWI maps are necessarily discontinuous and the classical 
Krylov-Bogolyubov argument cannot be applied. 
The first step to overcome this difficulty was done by B.~Pires (2016) \cite{Pi}, 
where the existence was proven under some (uncheckable) technical 
assumptions.\footnote{Basically 
       this is a version of the Krylov-Bogolyubov theorem, specialized 
       for discontinuous maps under study.} 
Let us mention also two other papers \cite{FD,Tr}, where the analysis again 
was based on a number of (uncheckable) technical assumptions. 
The first and the only real result about the existence of invariant measures 
of an orientation preserving circle PWI was obtained by S.~Kryzhevich (2020) \cite{Kr} 
by means of rational approximations. We will discuss this approach in some detail 
in Section~\ref{s:setup}.

\section{Setup}\label{s:setup}

Let $X$ be a subset of the Euclidean space $\IR^d,~d\ge1$ equipped with 
a certain metric $\rho(\cdot,\cdot)$ belonging to the class of 
$\t\rho_p(x,y):=(\sum_i|x_i - y_i|^p)^\frac1p,~p\ge1$, 
and let $\{X_i\}$ be a partition of $X$ into disjoint regions. 
By a {\em region} we mean a convex set with a nonempty interior.
The union of boundary points of all regions we denote by $\Gamma$, 
which by definition is of zero $d$-dimensional Lebesgue measure. 

\bdef{A {\em piecewise isometry (PWI)} is a map $T:X\to X$, satisfying the 
property that its restriction to each region $X_i$ is an isometry. 
We will refer to $\{X_i\}$ as a {\em special} partition, and to the 
maps $T_i|_{X_i}$ as {\em local} maps.}

\bdef{A restriction $T_{|X_i}$ is said to be an {\em extendable 
isometry} if it can be extended to an isometry $T_i$ on the entire 
$X$ such that $T_i|_{X_i}\equiv T_{|X_i}$. Correspondingly a 
piecewise isometry $T$ is said to be {\em extendable} if its 
restriction to each region $X_i$ is an extendable isometry.}

As we will see the innocent looking assumption of extendability is in fact 
somewhat restricting. Collecting well known geometric facts (see e.g. 
\cite{NDF, AKT}) we get the following

\begin{proposition}\label{p:geom} 
An isometry from a region $Y\subset \Tor^d$  to $\Tor^d$ is extendable 
if and only if it can be represented as a finite superposition of 
basic torus isometries: translations and coordinate flips or exchanges. 
\end{proposition}

The main result of our study is the following theorem about extendable 
torus piecewise isometries. 

\begin{theorem}\label{t:main}
Let $T$ be an extendable PWI of the unit torus $X$ with 
a finite special partition $\{X_i\}$. Then there exists at least one probabilistic 
$T$-invariant measure $\mu_T$. Additionally, if the boundary set $\Gamma$ 
contains no periodic points, this measure is non-atomic.
\end{theorem}

We expect that the extendability assumption can be avoided or significantly weaken, 
but the methods used in this article make extensive use of this assumption.
In the one-dimensional orientation preserving case the claim of Theorem~\ref{t:main} 
was proven in \cite{Kr}, but in the general (even one-dimensional) setting this 
claim is new. Moreover, in \cite{Kr} the non-atomic property of invariant measures 
has been proven without the assumption of the absence of periodic points at $\Gamma$. 
One of the key-points there was the demonstration that a periodic trajectory of a PWI cannot 
be isolated.\footnote{A trajectory is isolated if all nearby trajectories have other itineraries.} 
It was expected that even in the one-dimensional non orientation preserving case 
the situation will be the same. This turns out to be wrong, and a counterexample will 
be discussed in the last part of the proof of Proposition~\ref{p:local}.

To emphasize that the claim of Theorem~\ref{t:main} is far from being obvious,  
we will construct examples of extendable piecewise isometries without invariant 
measures, albeit with a countable special partition (see Section~\ref{s:absence}).

A general scheme of the proof of Theorem~\ref{t:main} is as follows. 
First, we consider a family of torus isometries $\{g_\theta\}$ satisfying the 
property that any finite subcollection of such isometries generates 
a group $G(g_{\theta_1},\dots,g_{\theta_n})$ having a finite number 
of elements. We refer to this property as weak periodicity.\footnote{See 
     Section~\ref{s:per} for discussion.} 
In Section~\ref{s:per} we will analyze properties of weakly periodic maps. 
In the case of commuting generators the situation is rather trivial and an 
important point here is that we are able to study some important cases 
when at least some of the the generators $g_\theta$ do not commute.

On the 2-nd step we consider torus piecewise isometries $T$ such that 
the restrictions $T_{X_i}$ can be extended to the torus isometries 
belonging to the above family. In Section~\ref{s:per} we study properties 
of such maps and prove for them the claim of Theorem~\ref{t:main}.

Finally in Section~\ref{s:appr} we approximate an arbitrary extendable torus 
piecewise isometry by means of weakly periodic ones considered in Section~\ref{s:per}, 
which will allow to finish the proof of Theorem~\ref{t:main}. 

It is worth noting, that this result does not contradict to the possibility 
of having attractors admitting singular invariant measures. 
Various numerical results (see e.g. \cite{Go,AKT}) give some evidence of 
the possibility of such behavior. 
The point is that the corresponding measures might not be represented as 
limits of invariant measures of maps discussed in Section~\ref{s:per} and hence 
cannot be controlled by our construction. 

A similar approach was used in \cite{Kr}, where a specialized version 
of the above construction was applied for the case of orientation preserving 
circle piecewise rotations. In that case the author used the approximation 
scheme based on rational translations and special partitions with rational 
boundary points. 
The last fact is a crucial point not allowing to use this approximation 
in the multidimensional setting. Indeed, for $d>1$ the boundaries of 
regions cannot consist of points with rational coordinates only.

\section{Focusing phenomenon and absence of invariant measures}\label{s:absence}
A natural approach to prove the existence of an invariant measure of a 
piecewise continuous dynamical system $(T,\Tor)$ is to consider a sequence of 
images $\{T^nm\}_n\ge0$ of the Lebesgue measure $m$ on the torus $\Tor$. 
Let $\mu$ be a limit point of this sequence. If the $\mu$-measure of a set $\Gamma$ 
of points of discontinuity of the map $T$ is 0, then under some simple technical 
assumptions the measure $\mu$ is $T$-invariant. The opposite situation, when 
the images of the Lebesgue measure concentrates on a small set we refer 
as the focusing phenomenon. To be precise let us introduce the {\em focusing set}:
$$ S(T) := \liminf_{n\to\infty}T^nX ,$$
where $\liminf_{n\to\infty}A_n:=\cup_n\cap_{k\ge n}A_k$ is the set of points 
belonging to all but finitely many sets $A_n$. Then the ``size'' of the set $S(T)$ 
describes the focusing quantitatively.

The focusing or concentration phenomenon is well known for systems with (locally) 
stable periodic orbits. The focusing occurs exactly to these periodic orbits, and the 
main reason for this is the (locally) contracting nature of the map $T$. 
In the situation with piecewise isometries we do not observe the contraction 
anywhere even locally. 

In this section we discuss how the focusing may occur for a PWI. 
In distinction to other parts of this paper we allow the special partition to be countable. 
The focusing possibilities for the case of the finite special partition will be discussed 
in the end of the Section.
The main idea under the following two examples is that under the action of the maps 
under study any measure converges weakly (focuses) to a measure $\mu$ supported 
by a very small wandering set. 
Recall that the latter means that all its forward images under dynamics do not 
intersect with the set. Therefore the measure $\mu$ cannot be invariant.  
Examples of this sort are known in the case of piecewise contractive maps, 
where the phenomenon of the measure concentration is natural (see e.g. \cite{Bl}), 
but look somewhat surprising for the piecewise isometry.

\begin{example}\label{ex:2} 
Let $X:=[0,1)^2$ and the special partition $\{X_n\}$ be defined as follows:
$X_n:=(1-2^{-n+1},1-2^{-n}]\times[0,1),~~n\ge1$, while 
$Tx:=(x_1+2^{-n-1},x_2)$ if $x\in X_n$. 
%The two-dimensional ($d=2$) example is demonstrated in the \bcr{Fig.~\ref{f:absence}}.
\end{example}

\begin{figure} \label{f:absence} \begin{center}
\begin{tikzpicture}[scale=0.4]
       \draw (0,0) rectangle (10,10); 
       \node at (0,-0.5) {$0$};  \node at (10,-0.5) {$1$};  \node at (-0.3,9.6) {$1$};
       \draw (10-5,0) rectangle (10-10,10); \node at (10-9,8) {$X_1$};
       \draw[fill=yellow] (10-2.5,0) rectangle (10-5,10);  \node at (10-4.5,8) {$X_2$};
       \draw (10-1.25,0) rectangle (10-2.5,10);        \node at (10-2,8) {$X_3$};
       \draw (10-.62,0) rectangle (10-1.25,10);    
       \draw (10-.31,0) rectangle (10-.62,10);
       \draw (10-.15,0) rectangle (10-.31,10);
       \node at (10-5,-0.7) {$\frac12$}; \node at (10-2.5,-0.7) {$\frac34$}; 
       \node at (10-1.25,-0.7) {$\frac78$};
       \bwait2{
       \draw[fill=green] (10-1.25,0) rectangle (10-3.75,10);  \node at (10-2.5,8) {$TX_2$};
       \draw (10-2.5,0) rectangle (10-5,10); }
      \bwait2{{\node at (17.7+1,8) {$X_n:=(1-2^{-n+1},1-2^{-n}]\times[0,1)$}; 
                          \node at (16.9+1,6.5) {$X:=\cup_{n\ge1} X_n = [0,1)\times[0,1)$}; 
                          \node at (17.2+1,5) {$Tx:=(x_1+2^{-n-1},x_2)$ if $x\in X_n$};
                          \node at (15.3+1,3.5) {$T^kX \toas{k\to\infty} \{0\}\times[0,1)$}; }}
\end{tikzpicture}
\end{center}
\caption{Example of a piecewise torus isometry without invariant measures. 
The element $X_2$ of the special partition and its image are marked in yellow and 
green colors correspondingly.}
\end{figure}
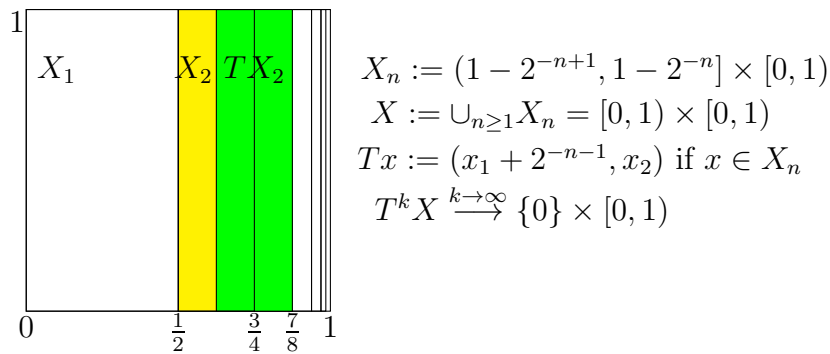

Observe that after exactly two iterations each element of the partition $\{X_n\}$ will be 
mapped to the union of elements with larger indices. 
In this situation obviously any measure converges under the action of $T$ to a measure 
supported by the unit segment $\{0\}\times[0,1)$. This segment is an wandering set, 
since its forward images under the action of $T$ do not intersect with the segment,  
which proves the absence of invariant measures. 

Nevertheless, in the Example~\ref{ex:2} we do not have a complete concentration of measure, 
namely  the situation when the limit (but not invariant) measure is a $\delta$-measure at 
a single point. To achieve this goal instead of the basically one-dimensional example above, 
we consider a bit more involved $d$-dimensional map. 

\begin{example}\label{ex:3}
Let $X:=[0,1)^d$ and the special partition $\{X_n\}$ be defined as follows:
$$ d(1-\frac1n) \le \sum_{i=1}^d x_i < d(1-\frac1{n+1}) ,$$
while $Tx:=x + \frac1{2n(n+1)} \bar1 \mod1$ for $x\in X_n$, 
where $\bar1$ is the $d$-dimensional vector with unit coordinates.
\end{example}

\begin{figure} \label{f:absence2} \begin{center}
\begin{tikzpicture}[scale=0.4]
       \draw (0,0) rectangle (10,10); 
           % \fill[green] (90:4) -- (210:4) -- (-30:4) -- cycle;
           % \fill[red] (0,0) -- (5,0) -- (0,5) -- cycle;
       \draw (0,0) -- (10,0) -- (0,10) -- cycle;  \node at (1,8) {$X_1$};
       \draw[fill=yellow] (0,10) -- (3.3,10) -- (10,3.3) -- (10,0) -- cycle; \node at (2,9) {$X_2$};
       \bwait2{
         \draw[fill=green]  (0+1,10+1) -- (3.3+1,10+1) -- (10+1,3.3+1) -- (10+1,0+1) -- cycle; 
         \node at (4.3,9) {$TX_2$};}
       \draw[thick] (5.5,10) -- (10,5.5); \draw[thick] (7,10) -- (10,7); 
       \draw[thick] (8.2,10) -- (10,8.2); \draw[thick] (9.1,10) -- (10,9.1);
       \node at (7.2,9) {$X_4$};
       \node at (0,-0.5) {$0$};  \node at (10,-0.5) {$1$};  \node at (-0.3,9.6) {$1$};
       \bwait2{{\node at (21.1,8) {$X_n:= \{x:~2(1-\frac1n) \le x_1+x_2 < 2(1-\frac1{n+1})\}$};
                          \node at (16.9+1,6.5) {$X:=\cup_{n\ge1} X_n = [0,1)\times[0,1)$}; 
                          \node at (19.3,5) {$Tx:=x + \frac1{2n(n+1)} \bar1 \mod1$ if $x\in X_n$};
                          \node at (15.3+1,3.5) {$T^kX \toas{k\to\infty} \{1\}\equiv \{0\}$}; }}
\end{tikzpicture}
\end{center}
\caption{Example of a piecewise torus isometry converging to $\bar1$. 
The element $X_2$ of the special partition and its image are marked in yellow and 
green colors correspondingly.}
\end{figure}

In distinction to the previous example the image of $X_n$ 
may intersect with elements $X_k$ with smaller indices $k<n$. Since the diameters of 
the elements of the partition decrease when their indices grow, after the implementation 
of the operation $\mod1$ peripheral parts of the image may be located  
below the central part and intersect with the elements with lower indices. 
%(see Fig.\ref{f:absence2}). 
In fact, $TX_n\cap X_1\ne\emptyset~\forall n$. 
Nevertheless the Hausdorff distance from the set $T^kX_n$ 
to the point $\bar1$ goes to 0 when $k\to\infty$. 
Therefore under the action of $T$ any measure weakly converges to $\delta_{\bar1}$, 
while the point $\bar1\equiv 0$ is again an wandering set.

\bigskip

If the special partition is finite we cannot provide examples with the focusing 
to the set of zero Lebesgue measure. Instead we try to estimate the ``size'' 
of the focusing set from bellow. 

\begin{example}\label{ex:f}
Let $X:=[0,1)$ and the special partition $\{X_n\}$ consists of $N$ intervals 
of equal length. Let $\alpha:=\{\alpha_i\}_{i=1}^n$ be a collection of numbers 
and let $T_\alpha x:=x + \alpha_i \mod1$ for $x\in X_i$. We are interested 
in the focusing for this family of maps, namely in 
$$ \kappa:=\inf_\alpha m(S(T_\alpha)) .$$
\end{example}

For a given $\alpha$, the set $T_\alpha$ is a typical one-dimensional orientation 
preserving piecewise isometry. Siince the images of the elements of the special 
partition may intersect, the value of $\kappa$ needs not to be close to one. 
However, since the map is locally isometric, by the same telescoping construction 
as in the examples~\ref{ex:2}, \ref{ex:3} one can achieve $\kappa=\frac1N$. 
From the first sight it seems that this estimate cannot be significantly improved. 
To demonstrate that this is absolutely not the case choose an arbitrary small 
parameter $\ep>0$ and set $\alpha_1:=\ep, \alpha_2:=-\ep, \alpha_n:=-\frac1N$ 
for $n>2$. In other words, $T_\alpha X_{n+1}=X_n$ for $n>3$ (telescoping), 
but the images of the first two elements intersect by the set of length $2\ep$. 
It is easy to calculate that $S(T_\alpha):=[\frac1N-\ep, \frac1N+\ep]$. 
Since $\ep>0$ is arbitrary, we get $\kappa=0$.

\section{Weakly periodic semigroups of maps}\label{s:per}

The main aim of this Section is to introduce the class of weakly periodic torus 
piecewise isometries and to prove Theorem~\ref{t:main} for them. 

\bdef{A collection of maps $\{g_i\}_{i=1}^n$ from a set $X$ into itself is said to be 
{\em weakly periodic} (notation WPer) if the number of elements in the 
semigroup generated by them is finite, i.e. $\#G(g_1,\dots,g_n)<\infty$.}

The reason under this terminology is that a periodic homeomorphism $T$ (i.e. $T^n=Id$) 
generates a semigroup with a finite number of elements, i.e. $T\in$ WPer.

The simplest (and very instructive) example of a weakly periodic semigroup of 
maps is a semigroup generated by a finite number of rational circle translations.

An important observation about this property gives the following result.

\begin{proposition}\label{p:nofim} 
Let $g_i\circ g_j=g_j\circ g_i~\forall i,j$, then 
$\sum_i\#G(g_i)<\infty$ implies $\#G(g_1,\dots,g_n)<\infty$. Otherwise 
there exists a pair of non-commuting maps $g_1,g_2$ such that $\sum_{i=1}^2\#G(g_i)=4$, 
but $\#G(g_1,g_2)=\infty$.
\end{proposition}

\proof Each element $T\in G(g_1,\dots,g_n)$ may be represented in the form 
of a finite superposition of the generators $g_i$, and due to their commutativity we get 
$$ Tx=g_1^{k_1}\circ \dots \circ g_n^{k_n} x,~~0\le k_i=k_i(x)<\infty ~~\forall x\in X .$$ 
By the assumption $\exists N<\infty$ such that $\#G(g_i)<N~\forall i$. 
Therefore, $$\#(\cup_{T\in G}Tx) \le N^n<\infty ~~\forall x\in X,$$
which implies the first claim. 

The 2nd claim is not so obvious as one might expect. In particular, in the one-dimensional 
case it is wrong for orientation preserving maps. The latter property turns out to be crucial 
here. Let $X:=[0,1),~g_1(x):=-x \mod1,~ g_2(x):=-x+\alpha \mod1,~\alpha\notin\IQ$. 
Observe that $\#G(g_i)=2, i=1,2$, but $(g_1g_2)^nx=x-n\alpha \mod1$, which implies
$\#G(g_1,g_2)=\infty$ since $\alpha\notin\IQ$. \qed

The definition of the weakly periodic maps implies the following simple technical result.

\begin{proposition}\label{p:ex-per} 
Let $\{g_i\}_{i=1}^n\in$ WPer and let a PWI map $T:X\to X$ with 
a finite special partition $\{X_i\}_{i=1}^n$ such that 
$T|_{X_i}\equiv g_i|_{X_i}~\forall i$. Then $\forall x \in X$ 
the trajectory $\{T^nx\}_{n\ge0}$ starting at $x$ is eventually periodic.
\end{proposition}

\proof The claim follows immediately from the observation that 
the total number of different points obtained by succesive applications 
of the generators $g_i$ is finite. \qed

Despite that according to Proposition~\ref{p:nofim} in general the WPer property 
for individual maps does not imply this property for a collection of them, 
the situation becomes more optimistic when we restrict ourselves to torus 
isometries. 

As we already noted in Proposition~\ref{p:geom}, any $d$-dimensional torus isometry 
may be represented by a superposition of basic isometries of the following 3 types: 
\begin{itemize}
\item[(i)] translation: $x\to x+ v \mod1,~~v\in\IR^d$
\item[(ii)] coordinate flip: $x_i \to -x_i$
\item[(iii)] coordinate exchange: $x_i\leftrightarrow x_j$
\end{itemize}

The maps belonging to different types do not commute, but at least the maps from 
the first two types commute with the maps of the same type, however the maps of type (iii) 
do not commute even between themselves. Weak periodicity of the maps from 
the first type is equivalent to their periodicity, while the maps of the other two 
types are always 2-periodic. 

To apply the machinery developed in Proposition~\ref{p:wper0} one needs to show 
that a PWI generated by a finite collection 
of periodic basic isometries satisfy the weak periodicity property. 

\begin{theorem}\label{t:wper} 
Let $T$ be a torus PWI such that $T|_{X_i}$ is a finite composition 
of periodic basic isometries ~$\forall i$. Then the map $T \in$ WPer.
\end{theorem}

The proof of this intuitively obvious claim (to which we devote the remaining part of 
this Section) turns out to be surprisingly cumbersome. 
In particular, we are not able to perform the proof directly for the map $T$ itself 
and are instead embed it to a certain semigroup generated by symmetrized torus 
extensions of the local isometries.

\bdef{By a {\em symmetrized semigroup} $G^*(g_1,\dots,g_n)$ generated by 
basic torus isometries $\{g_i\}$ we mean a semigroup generated by the maps 
$\{g_i\}$, where each translation map $x\to x+ v \mod1$ is considered 
together with all permutations and sign changes of the coordinates of the vector $v$, 
similarly each coordinate flip or exchange is considered together with all other 
coordinate flips or exchanges respectively. 
In the same manner by $(\{g_i\})^*$ we mean a symmetrized collection 
of basic torus isometries $\{g_i\}$.}

Fortunately to the finite dimension of the torus, the number of additional 
generators in the symmetrized semigroup is finite as well.\footnote{It cannot 
                                               be enlarged in more than $(d!)^3$ times.} 

\begin{proposition}\label{p:wper} Let $\{g_i\}_{i=1}^n$ be a collection 
of periodic basic torus isometries. 
Then $\#G^*(g_1,g_2,\dots,g_n) \le \prod_{i=1}^n \# G^*(g_i) <\infty$.
\end{proposition}

\proof We proceed by induction. The claim for $n=1$ is trivial, except 
for its finiteness. If $g_i$ belongs to (i) or (ii) types, then the finiteness follows 
from the periodicity assumption together with the commutation property. 
If $g_i$ is of the 1st type, the last argument cannot be applied, 
but instead one gets a direct estimate $G^*(g_i) \le d!$ since all possible 
coordinate permutations are taken into account. 

Let us make the induction step from $n$ to $n+1$. 
Observe that the application of the semigroup $G^*(g_1,g_2,\dots,g_n)$ to a given 
point $x\in X$, i.e. the set $G^*(g_1,g_2,\dots,g_n)x$ is invariant with 
respect to each of the maps $(\{g_i\}_{i=1}^n)^*$. The latter means that  
$$ g(G^*(g_1,g_2,\dots,g_n)x) = G^*(g_1,g_2,\dots,g_n)x \qquad  
                                                        \forall g\in (\{g_i\}_{i=1}^n)^* .$$ 
We will show that for any periodic basic isometry $g_{n+1}$ the set 
$G^*(g_{n+1})G^*(g_1,g_2,\dots,g_n)x$ is invariant with 
respect to each of the maps $(\{g_i\}_{i=1}^{n+1})^*$, which implies the claim. 

Denote $A:= G^*(g_1,g_2,\dots,g_n)x,~B:=G^*(g_{n+1})G(g_1,g_2,\dots,g_n)x$ and 
consider all possibilities:

\begin{itemize}

\item $g_{n+1}$ is a translation. If $A$ is invariant with respect to some 
translation $g_i$, then the the invariance of $B$ with respect to $g_i$ follows from the 
commutation property. The invariance of $A$ with respect to a coordinate flip is 
equivalent to the fact that the set $A$ is a uniform lattice on the torus. The latter 
property cannot be destroyed by a torus translation. 
If $A$ is invariant with respect to a coordinate exchange, then this set is 
characterized by the symmetry with respect to the corresponding unit cube diagonal. 
Unfortunately the translation $g_{n+1}$ {\em may destroy} this invariance for $B$. 
Indeed a symmetric pair $(a,b), (b,a)$ may by mapped to $(a+c,b), (b+c,a)$. 
This is exactly the point, where we need the symmetrization, because 
$G^*(g_{n+1})A$ clearly preserves the symmetry under question.

\item $g_{n+1}$ is a coordinate flip. If $A$ was invariant with respect to some 
translation, $B:=G^*(g_{n+1})A$ preserves this invariance due to the symmetrization. 
The invariance of $A$ with respect to a coordinate flip $g_i$ implies the 
same property for $B$, since coordinate flips commute. 
The invariance of $A$ with respect to a coordinate exchange is inherited by $B$ 
again helps to the symmetrization with respect to the coordinate flips.

\item $g_{n+1}$ is a coordinate exchange $x_i\leftrightarrow x_j$. This situation 
can be analyzed in exactly the same manner and we skip it.

\end{itemize}

The induction step is proven. \qed

This result finishes the proof of Theorem~\ref{t:wper}. Indeed, denote the torus 
extensions of the local isometries $T|_{X_i}$ by $g_i$ and assume that 
the maps $g_i$ are finite superpositions of periodic basic torus isometries. 
Then the trajectory of length $n$ under the action of $T$ of any given point $x\in X$ 
is a subset of $G^*(g_1,\dots,g_n)x$, which is uniformly finite on $x$ and $n$ 
by Proposition~\ref{p:wper}. This proves the weak periodicity of the map $T$. \qed

Additionally from the weak periodicity of the map $T$ in this setting we get 
some important information about its invariant measures.

\begin{proposition}\label{p:wper0}
Let $\{g_i\}_{i=1}^n$ be a collection of torus diffeomorphisms and 
let $\#G(g_1,\dots,g_n)<\infty$.
Let $T$ be a piecewise isometry of the unit torus $X$ with the 
finite partition $\{X_i\}$ and such that each $T|_{X_i}$ can be extended 
to a torus map $T_i\in G(g_1,\dots,g_n)$. 
Then there exists at least one absolutely continuous invariant measure, 
being a restriction of the Lebesgue measure to a certain region.
\end{proposition}

\proof Denote by $\t\Gamma:=\cup_{g\in G (g_1,\dots,g_n)}g^{-1}\Gamma$ -- 
the union of all pre-images under the action of all maps $g\in G(g_1,\dots,g_n)$ 
of the boundary $\Gamma$ of the elements of the partition $\{X_i\}$. 
By the weak periodicity assumption $\t\Gamma$ divides the torus $X$ into a 
finite number of regions $\{Y_j\}$ of positive Lebesgue measure, such that 
points belonging to the same region $\{Y_j\}$ have the same 
itineraries under the action of the map $T$. Therefore each of 
the regions $Y_i$ is either periodic under the action of $T$, or is wandering. 
The latter case is out of interest for us, while a collection of periodically 
exchanging regions $Y_{i_1},\dots,Y_{i_k}$ supports the Lebesgue measure 
on them being invariant with respect to the map $T$. Moreover, all absolutely 
continuous $T$-invariant measures may be described this way. \qed

\section{Approximation} \label{s:appr}

In what follows to control the rate of recurrence we will need the following 
abstract result. 

\begin{lemma}\label{l:kac} Let $(T,X,\cB,\mu)$ be a measurable dynamical 
system, let $A\in\cB,~\mu(A)>0$, and let $\tau_A(x)$ be the first return time 
for a point $x\in A$ to return to the set $A$ under the action of $T$. 
Then $\mu(\{x\in A:~ \tau_A(x)\le \frac1{\mu(A)}\})>0$.
\end{lemma}

Before to prove this result, observe that according to the famous Kac Lemma (see e.g. \cite{Si}) 
in the ergodic case the mathematical expectation of the first return time is 
exactly equal to $1/\mu(A)$. On the other hand, in our setting we cannot a priori 
assume ergodicity without which the claim of Lemma~\ref{l:kac} is as best as possible.  

\bigskip

\proof By the Poincare Recurrence Theorem $\mu$-a.a. points 
return to $A$ under the action of $T$. 
Let $A_n:=\{x\in A:~n\le \tau_A(x)<\infty\}$ be the set of points returning 
to $A$ at least after $n\ge0$ iterations, and suppose 
$A\setminus A_n=\emptyset \hence \mu(A)=\mu(A_n)$.  
Since $T^{-k}A_n \cap A =\emptyset~~\forall k<n$, we get  
$T^{-i}A_n \cap T^{-j}A_n = \emptyset~~\forall i\ne j<n$ 
(otherwise the points from the intersection above will not return to A). 
Hence 
$$ 1 \ge \mu(\cup_{k=0}^{n-1}T^{-k}A_n) = \sum_{k=0}^{n-1}\mu(T^{-k}A_n) 
                                                                = n\mu(A) $$
and thus $n\le 1/\mu(A)$. \qed

Let $\{T^{(n)}\}$ be a sequence of piecewise torus isometries with the common 
finite special partition $\{X_i\}, ~\Gamma:=\cup_i\partial X_i$, 
satisfying the assumptions of Proposition~\ref{p:wper0}, 
and let $\{\mu_n\}$ be a sequence of their absolutely continuous invariant measures.

\begin{lemma}\label{l:apr} Let $z\not\in\Per(T)$. Then 
$\forall\ep>0~\exists\delta>0 ~\limsup\limits_{n\to\infty}\mu_n(B_\delta(z))<\ep$.
\end{lemma}

\proof Assume from the contrary that there is $\ep>0$ and a pair of sequences 
$\{n_k\}\toas{k\to\infty}\infty$ and $\{\delta_k\}\toas{k\to\infty}0$ such that 
$\mu_{n_k}(B_{\delta_k}(z))\ge\ep$. 
Then by Lemma~\ref{l:kac}
$$ \exists N<1/\ep:~~ (T^{(n_k)})^NB_{\delta_k}(z) \cap B_{\delta_k}(z) \ne \emptyset ,$$
which implies that $z\in {\rm Per}(T)$. We came to the contradiction. \qed 

Since we deal with a compact phase space, the sequence of the above mentioned 
absolutely continuous invariant measures of the approximating weakly periodic maps 
has a limit point. Denote this limit point by $\t\mu$ and choose a sequence of the 
measures  $\{\mu_n\}$ converging to it in the weak topology.

\begin{lemma}\label{l:noPer} Let $\Gamma\cap {\rm Per}(T)=\emptyset$. 
Then $\t\mu(\Gamma)=0$ and $\forall k, \phi\in C^0(X,\IR)$ we have 
$$ \lim_{n\to\infty}\mu_n(1_{\Clos(X_k)}\cdot\phi) = \t\mu(1_{\Clos(X_k)}\cdot\phi) .$$ 
\end{lemma}
\proof The first statement follows from Lemma~\ref{l:apr}. 
To prove the second statement without loss of generality we assume that $|\phi(x)|\le1$. 
Using again Lemma~\ref{l:apr} we get: $\forall \ep>0 ~\exists \delta>0$ such that 
$\mu_n(B_\delta(\partial X_k))<\ep$. Choose now a nonnegative function $h\in C^0(X,[0,1])$ 
such that 
$h_{|X_k\setminus B_\delta(\partial X_k)}\equiv1$ and 
$h_{|X \setminus B_\delta(\partial X_k)}\equiv0$. 
Then $\lim_{n\to\infty}\mu_n(h\phi) = \t\mu(h\phi)$. 
by the definition of the weak convergence of measures.

Therefore 
$$ \lim_{n\to\infty}|\mu_n(\phi) - \t\mu(\phi)| \le m(B_\delta(\Gamma)) ,$$ 
which implies the result since the value of $\ep>0$ can be taken arbitrarily small. \qed

\begin{lemma}\label{l:noPer1} Let $\Gamma\cap {\rm Per}(T)=\emptyset$ and 
let $\{T^{(n)}\}$ be a sequence of approximating weakly periodic maps, whose 
invariant measures $\mu_n\toas{n\to\infty}\t\mu$. Then
$$\mu_n(\phi\circ T^{(n)} - \phi\circ T)\toas{n\to\infty}0~~\forall \phi\in C^0(X,\IR) .$$
\end{lemma}
\proof Given $\ep>0$ by Lemma~\ref{l:apr} we can choose $\delta>0$ such that 
$\mu_n(B_\delta(\Gamma))<\ep$.\footnote{Note that the convexity of the elements 
      of the special partition is important at this point. Otherwise one needs to assume 
      that the boundaries of those elements are smooth enough.}
From the approximation assumption $\exists n_0$ such that 
$|\phi\circ T^{(n)}(x) - \phi\circ T(x)|<0$  $\forall n>n_0$ 
and $x\not\in B_\delta(\Gamma)$. This implies the claim. %\bcr{Check?} 
\qed 

Now we are ready to finish the proof of Theorem~\ref{t:main}. 

First, observe that each map, satisfying the assumption of Theorem~\ref{t:main}, 
can be approximated arbitrary well by weakly periodic ones. Indeed, by since 
$T|_{X_i}$ is extendable it can be represented by a finite superposition basic 
isometries. If there are isometries of type (i) among them, they can be approximated 
arbitrary well by rational translations (which are periodic). Basic isometries of 
other types are periodic from the very beginning.

The second step is to check that ``good'' invariant measures of the weakly 
periodic approximations, existing by Theorem~\ref{t:wper} and Proposition~\ref{p:wper0}, 
converge to a $T$-invariant measure. 

\begin{lemma}\label{l:fin} If $\exists \gamma\in \Gamma\cap {\rm Per}(T)$, 
then the measure uniformly distributed on the trajectory of $\gamma$ is $T$-invariant. 
Otherwise, if $\Gamma\cap {\rm Per}(T)=\emptyset$ the measure $\t\mu$ 
constructed in Lemmas above is $T$-invariant. 
\end{lemma}
\proof 
$$ |\t\mu(\phi - \phi\circ T)|  \le |\t\mu(\phi - \phi\circ T) - \mu_n(\phi - \phi\circ T)| $$
$$ \hskip4.7cm  +~ |\mu_n(\phi - \phi\circ T^{(n)})| + |\mu_n(\phi\circ T^{(n)} - \phi\circ T)| .$$
The first term in the rhs goes to 0 as $n\to\infty$ by Lemma~\ref{l:noPer}, 
the second term vanishes for each $n$ since the measures $\mu_n$ are $T^{(n)}$-invariant, 
while the third term goes to 0 by Lemma~\ref{l:noPer1}. 
Hence $\t\mu(\phi) = \t\mu(\phi\circ T)$, which implies the claim. \qed

\section{Torus restriction of a plane isometry} \label{s:affine}
For the family of torus metrics $\rho_p(\cdot,\cdot)$ for all $p\ne2$ the list 
of extendable isometries (translations, coordinate flips and coordinate exchanges) 
coincides with all possible isometries. However in the important Euclidean case 
$\rho_2(\cdot,\cdot)$ there is an additional isometry -- rotation, which is not 
extendable. The techniques developed in this paper does not allow to consider 
local rotations, because we do not have their weakly periodic approximations. 

In this section we consider actions of plane isometries restricted to the unit 
two-dimensional torus $X:=\Tor^2$ as the simplest models local rotations. 
Some results in this direction (mainly about topological properties) may be 
found in \cite{LV} (see further references therein). 
Recall that a general plane isometry can be written as $x\to \R x+b$, 
where $x,b\in \IR^2$ and $\R$ is an orthogonal matrix. Restricting this map 
to the torus we get the torus map $Tx:=\R x+b \mod1 : \Tor^2\to\Tor^2$.

\begin{example}\label{ex:1}
Let $\R=\R(\phi)$ be the rotation around the origin by the angle $\phi\in[0,2\pi)$ 
and let $b=0$. 
\end{example}

Having in mind that that each trajectory of the plane rotation belongs 
to a single circle around the origin, while a typical torus translation 
$x \to x+ b\mod1$ fills the torus densely, one naively expects something 
of the same sort for the map $T$. 
It turns out that if the rotation angle $\phi\ne0$ the 
situation is completely different (see Fig.~\ref{f:linear}). 

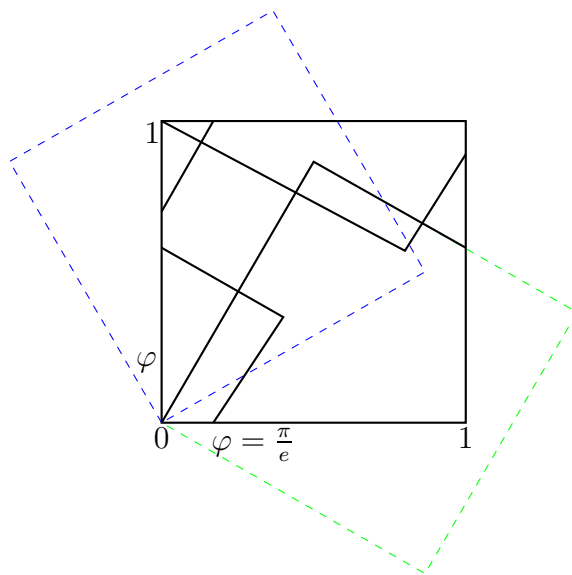
\begin{figure} \label{f:discontinuities} 
\begin{center}
\begin{tikzpicture}[scale=0.4]
       \draw[thick] (0,0) rectangle (10,10); 
       \node at (0,-0.5) {$0$};  \node at (10,-0.5) {$1$};  \node at (-0.3,9.6) {$1$};
       \bwait2{
       \draw[blue,dashed,rotate=+30] (0,0) rectangle (10,10); \node at (-.5,2) {$\phi$};} %dotted
       \bwait3{
       \draw[green,dashed,rotate=-30] (0,0) rectangle (10,10); \node at (3,-.72) {$\phi=\frac\pi{e}$};}
       \bwait4{
       \draw[thick] (0,0)--(5,8.65)--(10,5.8);
       \draw[thick] (0,10)--(8,5.7)--(10,8.9);
       \draw[thick] (0,5.8)--(4,3.5)--(1.7,0);
       \draw[thick] (0,7)--(1.7,10);}
\end{tikzpicture}
\end{center}
\caption{Thin blue/green lines -- the boundaries of forward/backward images of the boundaries 
of the unit square under the rotation by the angle $\pi/e$ around the origin. 
Thick black lines -- boundaries of the special partition.}
\end{figure}

The map $T$ defined above, considered as a map from the unit torus $X:=\Tor^2$ 
into itself is continuous everywhere inside the open unit square $\t{X}:=(0,1)^2$ 
and the set of discontinuities is represented by the boundary of this set $\partial\t{X}$. 
On the other hand, the same map may be considered as a map from the unit square 
$\hat{X}:=[0,1]^2$ into itself. In the latter case the set of discontinuous is much more 
involved (see Fig.~\ref{f:discontinuities}).

If $\det(I-\R)\ne0$ there is at least one fixed point of the map $Tx:=\R x+b \mod1$. 
Here and in the sequel by $I$ we denote the unit matrix. 
Indeed, making the change of variables $x=y+a$, we rewrite the equation for the 
fixed point of the map $T$ as $y+a=\R(y+a)+b \mod1$, which gives a solution 
$y:=0, ~ a:=(I-\R)^{-1}b \mod1$. Thus $z:=(I-\R)^{-1}b \mod1$ is the fixed point. 
If additionally, $b=0$ (case of a pure rotation) the origin is always the fixed point 
even without the assumption about the absence of the degeneracy of the matrix $I-\R$ 

Therefore the main question is about the existence of nontrivial periodic trajectories. 
Our numerical observations allow to formulate the following

\bigskip

\n{\bf Hypothesis.}\label{t:linear} {\em 
Let $Tx:=j\R(\phi)x + b \mod1$, where $j=\pm1,~\phi\in[0,2\pi), x,b\in\IR^2$. 
Then $\forall j,\phi,b~~\exists N=N(j,\phi,b)\in\IZ_+$ and $z=z(j,\phi,b,N)$ such that 
$T^Nz=z$. Moreover for a given triple $j,\phi,b$ one may find arbitrary large $N$ 
satisfying the above property.}

\bigskip

\n{\bf An idea of the possible proof.} 
We consider only the case of the pure rotation (i.e. $j=1,~b=0$) 
by the irrational angle (i.e. $\phi=\gamma\pi,~\gamma\notin\IQ)$,
leaving the general case (where a number of cumbersome combinatorial type 
calculations are necessary) for the future. 

Observe that $Tx = \R x \mod1 = \R x + \xi(x)$, where $\xi(x) \in\IZ^2$. Therefore 
to find a periodic point of period $1$ one needs to solve the equation:
$$ z = \R z + \xi_1 ,$$
and for the period $2$:
$$ z = \R (\R z + \xi_1) + \xi_2 = \R^2z + \R \xi_1 + \xi_2 $$
Continuing this, for a point of period $n$ we have  
$$ z = \R^nz + \sum_{j=1}^{n}\R^{n-j}\xi_j ,$$
which is equivalent to 
$$ (I-\R^n)z = \sum_{j=1}^{n}\R^{n-j}\xi_j.$$
If the angle $\phi$ is irrational, then the inverse matrix $(I-\R^n)^{-1}$ is 
well defined. Thus we have an explicit solution 
$$ z:=(I-\R^n)^{-1}\sum_{k=1}^{n}\R^{n-k}\xi_k .$$ 

Unfortunately this is not the end of the story, since in fact the integer vectors $\xi_i$ 
depend on $z$ and should be chosen in such a way that 
$\R^kz + \sum_{j=1}^{n}\R^{n-j}\xi_j$ 
belong to the unit square for each $1\le k\le n$, otherwise further corrections 
need to be taken into account. 

\bigskip

Let us formulate and prove a partial result leaving the general case for future studies.

\begin{proposition} Let $\phi\notin[-\phi_0,\phi_0]$, where $\tan\frac{\phi_0}2 = \frac12$. 
Then the equation $A(\phi) x = x$ always has a nontrivial solution.
\end{proposition}

\proof Consider a point $x$ lying on a vertical line with the 1st coordinate being equal to $\frac12$, 
i.e. $x=(\frac12,x_2)$. Then under the rotation around the origin by the angle 
$\phi:=-2\arctan(\frac1{2x_2})$ we get a point $\t{x}:=(-\frac12,x_2)$, which coincides  
with $x=(\frac12,x_2)$ modulo 1. The assumption $\phi\notin[-\phi_0,\phi_0]$ is explained
by the choice of the point $x$. 
\qed

\begin{proposition}\label{p:per} 
Let $z$ be a $N$-periodic point and let $R:=\inf_k\dist(T^kz,\partial X)>0$.
Then $S(r):=\cup_{i=0}^{N-1}B_r(z)$ is a forward invariant set $\forall r\in(0,R)$. 
\end{proposition}

\proof Consider any point $x$ close enough to $z$. Under the action of the map $T$ the 
$\rho_2(\cdot,\cdot)$ distance between $T^nx$ and $T^nz$ is preserved. Therefore each 
circle centered at $z$ with a small radius $r$ is mapped into itself. The value of the 
radius $r$ is bounded from above by the condition that the rotated circle does not 
intersect with the boundary of the unit square. If additionally $\phi/\pi$ is 
an irrational number, then the trajectory of the point $x$ fills densely the corresponding 
circle, otherwise, it is periodic. \qed

Observe that the number of ergodic invariant measures is always countable 
in distinction to the one-dimensional case, where it is known that in a ``typical'' situation  
the number of ergodic measures is bounded from above by the number of elements 
of the special partition.

A substantially more general result about the dynamics of a general PWI in a neighborhood of 
a given trajectory based on the Proposition~\ref{p:per} may be formulated as follows.

Let $T$ be a PWI of $X\in\IR^d$ with respect to a metric $\rho(\cdot,\cdot)$ and a special 
partition $\{X_i\}$.

\begin{proposition}\label{p:local} 
Let a trajectory of a point $z\in X$ be separated fromcthe boundary of 
the special partition $\Gamma$, i.e. $\inf_{n\ge0}\rho(T^nz,\Gamma):=R>0$. 
Then $\forall x\in B_R(z)$ the itineraries of the trajectories started at the points $x,z$ 
coincide, i.e. $T^nx\in X_i \Longleftrightarrow T^nz\in X_i ~\forall n\ge0$. 
Moreover, $\forall r\in(0,R)$ the spheres $S_r(T^n z)$ of radius $r$ centered at $T^nz$ 
are mapped one to another under dynamics, i.e. $TS_r(T^n z)=S_r(T^{n+1} z)$. 
If the separation property breaks down the trajectory stated at $z$ may be completely 
isolated in the sense that $\lim\inf_{n\to\infty}\rho(T^nx, T^nz)>0$ whenever $x\ne z$. 
\end{proposition}

\proof The separation assumption means that the ball $B_R(z)$ at each iteration of 
the map $T$ falls entirely to a single element of the special partition, which implies 
the claim about the itineraries. Now using the same argument as in the proof of 
Proposition~\ref{p:per} we see that the sequence of balls $B_R(T^nz)$ splits into 
the continuous on-parametric family of concentric spheres $S_r(T^n z)$ mapped 
one to another under dynamics. 

Observe that in the case of a periodic trajectory the separation property becomes 
especially simple and trivially checkable: $\cup_{n\ge0}T^n z \setminus \Gamma = \emptyset$. 

To demonstrate the last claim we return to the pure rotation of the 2-dimensional 
torus around the origin by the angle $\pi/4$. Clearly, the point $z=0$ is fixed under 
dynamics. Observe that a sufficiently small neighborhood $O_z$ of this point 
on the torus consists of 4 triangular sectors at the corners of the unit square and 
that each of these sectors belong to at most 2 elements of the special partition. 
This allows to use a simple numerical simulations to check that 
$T^n(X\setminus\{z\})\cap O_z=\emptyset~~\forall n\gg1$, which finishes the proof.  

Our numerical simulations indicate that the last claim remains valid for any rotation 
angle $\phi\ne0$ and any shift $b$, but this is not proven yet. Instead we consider 
an even more unexpected situation when a periodic trajectory of an one-dimensional 
non orientation preserving PWI is isolated. 

\begin{example}\label{ex:0}
Consider a PWI $T:X\to X:=[0,1)$ with the special partition into 4 intervals 
$X_1:=[0,\frac14), X_2:=[\frac14,\frac12),  X_3:=[\frac12,\frac34], X_2:=(\frac34,1)$,   
where the map is defined as follows $TX_1:=X_4, TX_2:=-X_3, TX_3:=X_1, TX_4=X_4$. 
The minus in $TX_2:=-X_3$ means the flip (orientation change) of the of the interval $X_3$, 
i.e. $T\frac14=\frac34, T\frac12=\frac14$. %(see Fig.\ref{f:isolation}). 
\end{example}

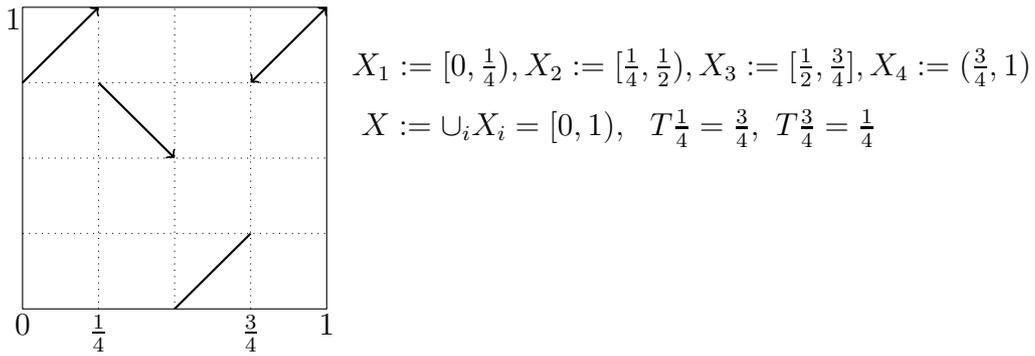
\begin{figure}
\begin{center}
\begin{tikzpicture}[scale=0.4]
       \draw (0,0) rectangle (10,10); 
       \draw[dotted] (2.5,0) -- (2.5,10); \draw[dotted] (5,0) -- (5,10); \draw[dotted] (7.5,0) -- (7.5,10); 
       \draw[dotted] (0,2.5) -- (10,2.5); \draw[dotted] (0,5) -- (10,5); \draw[dotted] (0,7.5) -- (10,7.5); 
       \draw[thick][->] (0,7.5) -- (2.5,10); \draw[thick][->] (2.5,7.5) -- (5,5); 
       \draw[thick][] (5,0) -- (7.5,2.5);   \draw[thick][<->] (7.5,7.5) -- (10,10);
       \node at (0,-0.5) {$0$};  \node at (10,-0.5) {$1$};  \node at (-0.3,9.6) {$1$};
       \node at (2.5,-0.8) {$\frac14$}; \node at (7.5,-0.8) {$\frac34$};
       \bwait2{{\node at (22,8) {$X_1:=[0,\frac14), X_2:=[\frac14,\frac12), 
                                                     X_3:=[\frac12,\frac34], X_4:=(\frac34,1)$};
                          \node at (19.6,6) {$X:=\cup_{i} X_i = [0,1), ~~T\frac14=\frac34,~T\frac34=\frac14$}; 
                   }}
\end{tikzpicture}
\end{center} \label{f:isolation} 
\caption{A non orientation preserving circle PWI with an isolated periodic trajectory $1/2, 3/4$.}
\end{figure}

Observe that in this example the trajectory started at the point $\frac14$ is 2-periodic,  
while all points from the surrounding intervals $(0,\frac14)$ and $(\frac14,\frac34)$ in 
2 iterations are mapped to the 4-th interval $(\frac34,1)$ and never return back. 
\qed

%%%%%%%%%%%%%%%%

\end{document}